# Piecewise closed form expression of Rössler-like trajectories

Stefano Morosetti
Department of Chemistry, University of Rome 'La Sapienza', Piazzale A. Moro 5, 00185 Roma, Italy
E-mail address: stefano.morosetti@uniroma1.it
http://orcid.org/0000-0002-3937-3045

*Abstract*

*Background*: Many dynamics are described by a set of nonlinear differential equations, and therefore can give rise to periodic and chaotic trajectories. Simple models such as Rössler dynamics have been introduced to study their behavior. Trajectories are obtained by numerical integration, or as infinite series solutions via Homotopy Analysis Method (HAM).

*Goal*: Get Rössler-like trajectories, by combining a carrier function with analytic piecewise functions. The iteration of the carrier function produces periodic and chaotic regimes as its parameter changes. The analytic piecewise functions describe the trajectory in space, using the values obtained from the carrier function.

*Method*: The first-return map of the maximum of the $X$ variable in the Rössler flow suggests the carrier function to be used. The study of trajectories of the numerically solved Rössler flow suggests analytic piecewise functions able to express the spatial coordinates over time.

*Results*: Rössler-like trajectories are obtained for periodic and chaotic regimes.

*Conclusions:* The use of iteration and piecewise functions allows analytic expression of the trajectories of an Rössler-like attractor, avoiding infinite series solution. It seems possible to extend this approach to other attractors, even if the differential equations are not known, but the experimental data have been collected.

*Introduction*

It is very common that physical, chemical or other phenomena from other fields of science are described by a set of nonlinear differential equations. As a result they can give rise to periodic and chaotic dynamics. In these cases the trajectories of the phenomenon under study, are not obtainable in analytic form according to Poincaré's theorem[1]. This means that the Hamiltonian of the system is not expressible by a series-developed perturbation of a known unperturbed Hamiltonian. However, a solution can be found in another way, as demonstrated by the three-body problem, not solvable according to Poincaré, but solved by Sundman[2] with a uniformly convergent series. More recently, the Homotopy Analysis Method[3] (HAM) has been developed to address these systems. In all cases a series solution must be calculated for each combination of parameters under study and necessarily truncated for practical reasons. This truncation can generate uncertainty about the correctness of the computed trajectories in the chaotic regime, due to their divergence resulting from positive Lyapunoff exponent. Alongside these analytical methods, trajectories are often obtained by numerical integration methods of the system of differential equations.

To describe the general properties of these systems, simple models such as those of Rössler[4] and Lorenz[5] are often used.

This article focuses on the Rössler model.
The Rössler flow is described by the following system of equations:

$$\dot{X} = -Y - Z$$
$$\dot{Y} = X + aY$$
$$\dot{Z} = b + XZ - cZ$$

where parameters $a$ and $b$ are commonly given the value 0.2 and $c$ is the control parameter that determines the flow regime. In fact, changing parameter $c$ passes through various regimes: successive duplications of the trajectory up to the chaotic regime. As an example, Fig. 1 shows the trajectories for c=4.00 (period 4 regime) and c=5.70 (chaotic regime).

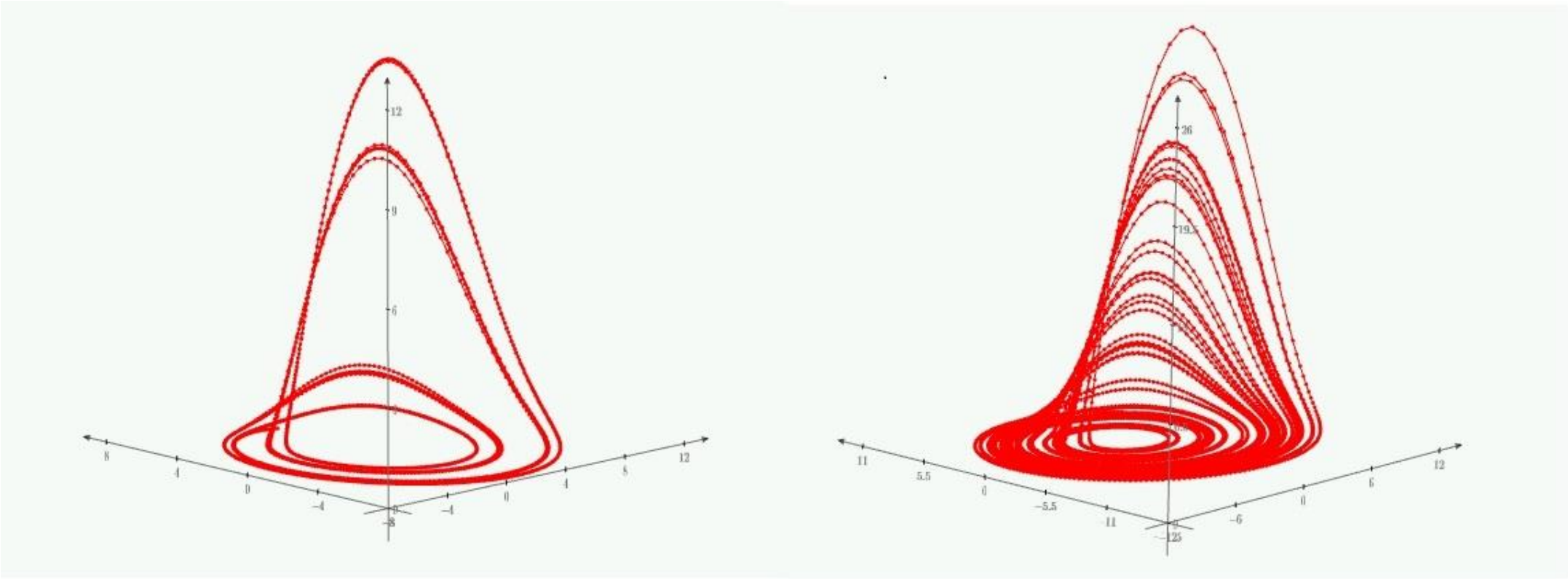

Fig. 1. Trajectories for $c$=4.00 (left) and $c$=5.70 (right) obtained through the PTC Mathcad Prime 3.0 program, using the AdamsBDF function to obtain numerical integration of the Rössler set of differential equations. The Mathcad Prime 3.0 worksheet used for the numerical integration of case $c$=4.00 is listed as S1 in *Mathcad Prime 3.0 worksheets* chapter.

This article proposes a mathematical apparatus capable of producing trajectories similar to those of Rössler, and in which there is a parameter that allows the transition from one periodic regime to another, until it reaches the chaotic one. These trajectories are not the solution of the Rössler flow, however, they represent the solution of a flow with Rössler-like and nonlinear behavior, as periodic and chaotic regimes are present.

*Method*

A function is used to get the transition from one periodic regime to another, until it reaches the chaotic one, as one parameter changes. It is called the carrier function. This also provides partial trajectory data.
Piecewise functions (one for each variable) are then obtained to fully describe the trajectories. They provide the coordinates. The expression of these functions is suggested by the study of the trajectories of the numerically solved Rössler attractor.
Carrier function and coordinates are described below.

<u>Carrier function</u>

From the time sequence of the values of maximum $X$, a graph can be obtained, showing each maximum $X$ value as a function of the previous one. It is called the first-return map of maximum $X$ and allows to derive the temporal sequence of these values (see Fig.2).

Changing the parameter $c$ of the Rössler flow changes the graph and thus passes through the various regimes.
Since the trajectories projected on plane $X, Y$ are approximately circular and centered around 0 (see Fig. 1), we assume that maximum $X$ also represents the radius $R$ of the trajectory, projected onto the plane $X$,$Y$, when it passes for the maximum $X$. In this process of approximate reconstruction of trajectories, we indicate a generic term for this sequence of radius values with $r_i$. There will be a $r_i$ for each subsequent revolved orbit around the origin of the trajectory projected onto the $X$,$Y$ plane.

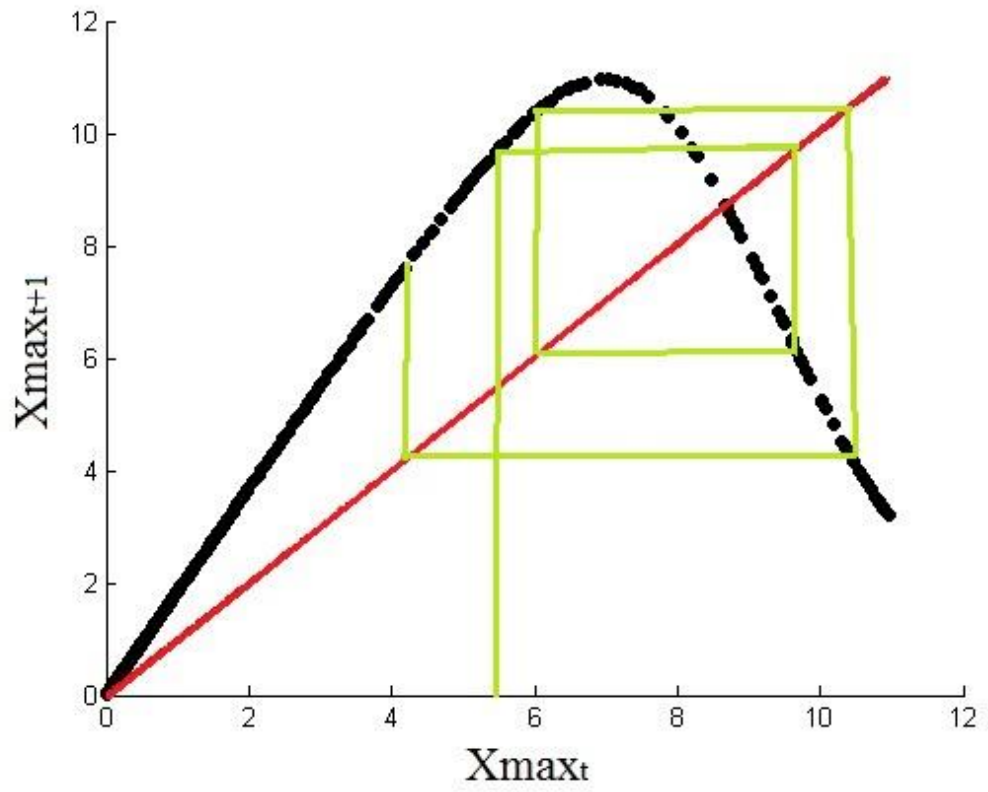

Fig. 2. In black the first-return map of maximum $X$. In red the diagonal line, used as a support to graphically find the next values. In green a possible sequence of successive values of maximum $X$, starting on the axis of $X_{max_t}$.

Note that the graph of the first return of maximum $X$ has a shape very close to the well-known curve, called the logistic map, $y = \lambda x(1 - x)$, apart from the scale. We will use this expression to obtain a sequence of values to which we will give the meaning $r_i$, taking into account that in this article we are interested in the reproduction of the trend of trajectories and not their exact numerical reproduction. It is known that the iteration of this function produces a sequence of numbers with a periodicity that depends on the parameter $\lambda$ and that a chaotic sequence can also be reached. Thus in iteration the parameter $\lambda$ assumes the same role as $c$ in the Rössler flow. Figure 3 shows the iteration for $\lambda = 3.5$ which produces a stable sequence of period 4.

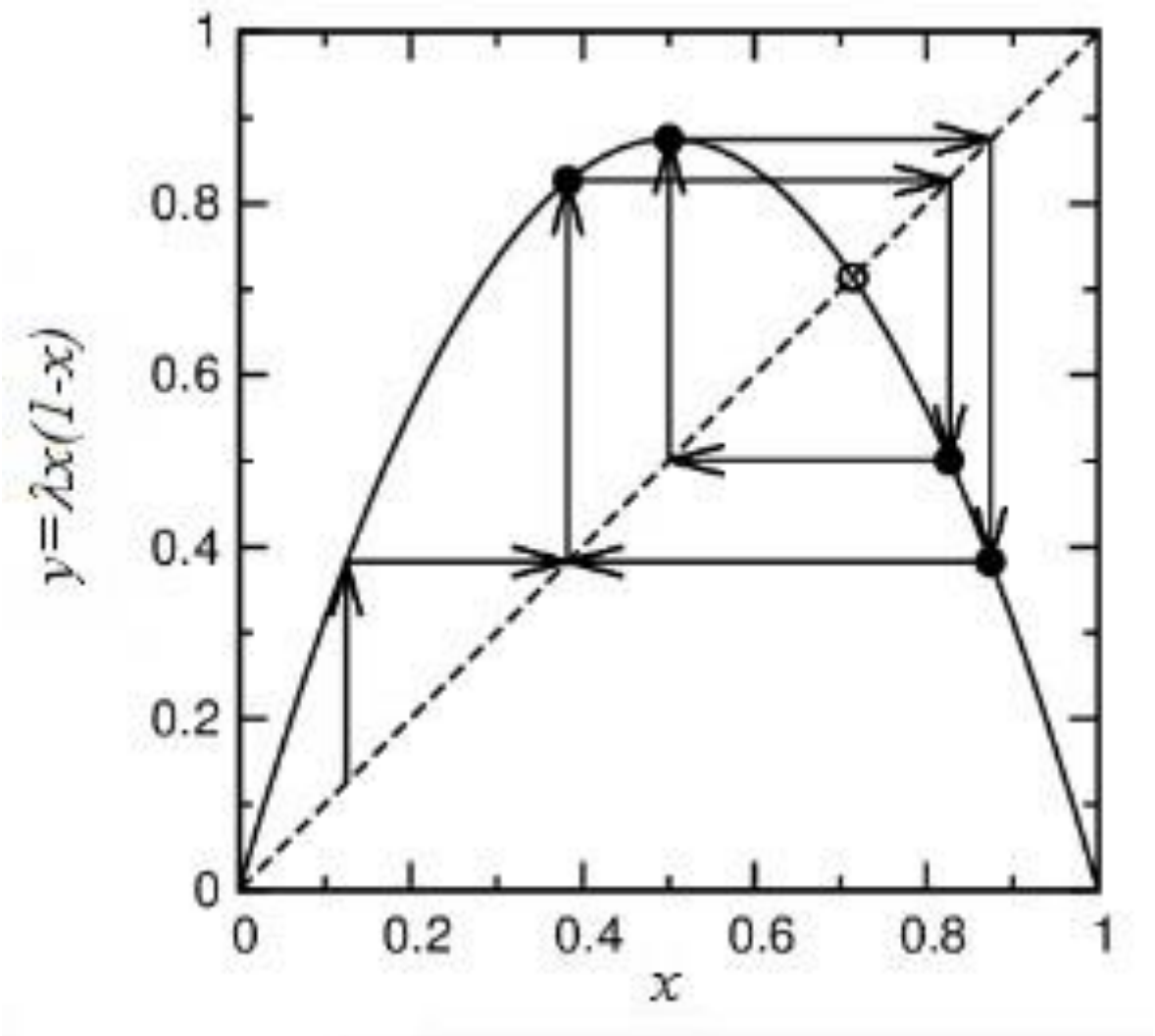

Fig. 3. Iteration of the function $y = \lambda x(1 - x)$, with $\lambda = 3.5$, which generates a period-4 sequence.

So from the iteration of this expression we get a sequence of values of $x$, to which we attribute the meaning of successive maximum $X$ values, by analogy with the Rössler's first return map. From what has been said, we attribute to these values the meaning of radius $R$ of the trajectory projected onto plane $X, Y$, in the moment it passes for the maximum $X$, thus obtaining a sequence of values $r_i$. Comparing the trend of the maximum $X$ values in the Rössler curves and the $x$-values obtained in the iterations of the parabola, the same pattern is found, that is, the same trend in relative values, even if the absolute values are different. This confirms the validity of the procedure. In Fig.4 the case of period-4 is presented.

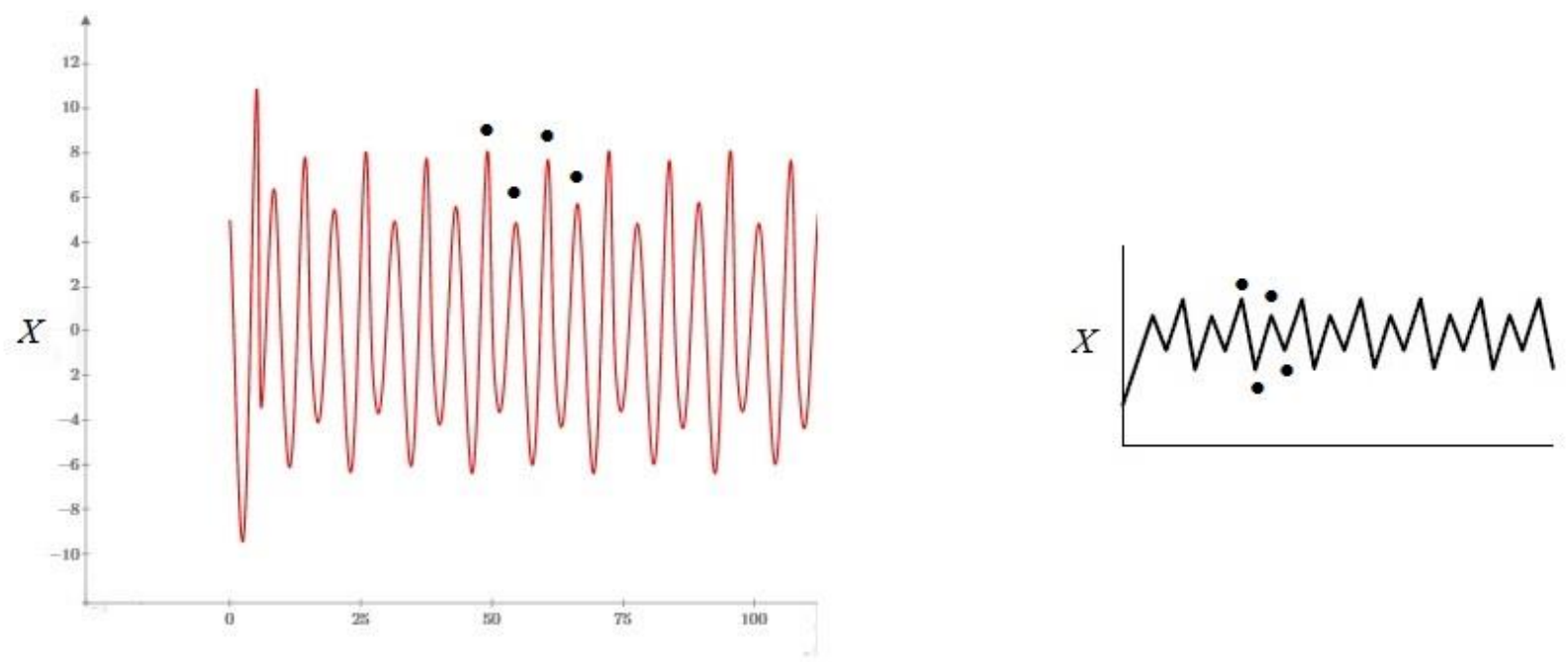

Fig.4: sequence of the $X$ values in the case of period-4. On the left the Rössler flow by $c = 4.00$. On the right, consecutive $X$ values, obtained by iteration of the function $y = \lambda x(1 - x)$ for $\lambda = 3.5$, are joined. In both cases, the values that follow one another in a cycle are highlighted with a point (only the maximum $X$ values are highlighted in the case of the Rössler flow). If 0, 1, 2, 3 is the increasing value of $X$, in both cases the temporal sequence is 3, 0, 2, 1.

The iteration of this expression can therefore act as a carrier function that allows to pass through cycles of different periodicity as the parameter $\lambda$ changes, and provides a sequence of values $r_i$ for each cycle (see eq. 1).

$$\begin{cases} x_i = \lambda x_{i-1}(1 - x_{i-1}) \\ r_i = 10 x_i \end{cases} \qquad \text{eq. 1}$$

The values $x$ are multiplied by 10 to have the same order of magnitude as the attractor's radii.

## Coordinates

The carrier function only obtains partial information on the trajectory. Reconstructing the trajectory requires: a) an expression that connects the values $r_i$ obtaining $R(t)$ on plane $X, Y$; b) to express $X$ and $Y$ as a function of $R(t)$; c) to find an expression for $Z(t)$.

To derive expressions that behave similar to Rössler's attractor, we start from its study. Trajectory is obtained by numerical integration of the attractor. $R(t)$, $Z(t)$, $X(t)/R(t)$ and $Y(t)/R(t)$ are reported, depending on the angle θ described on plane $X, Y$ from a point that proceeds over time along the trajectory, projected onto plane $X, Y$. Figure 5 shows the case for $c = 4.00$, but the other cases have similar graphs.

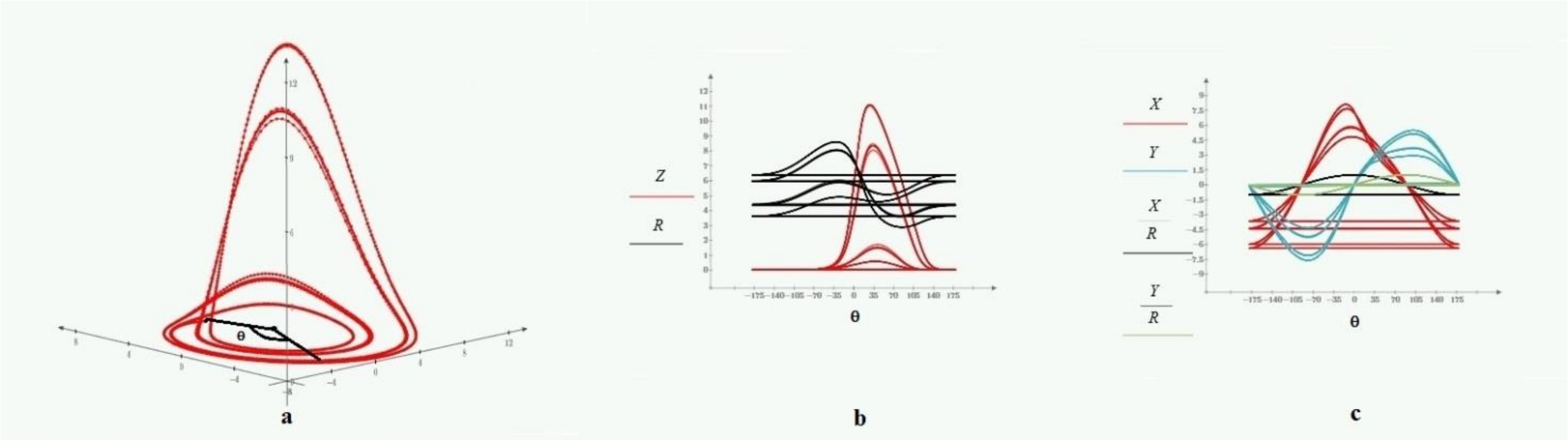

Fig.5: **a**) shows the angle θ representing the progression of the trajectory on plane $X, Y$. **b**) radius $R(t)$ on plane $X, Y$ and value $Z(t)$ as a function of θ. **c**) $X(t)$, $Y(t)$, $X(t)/R(t)$, $Y(t)/R(t)$ as a function of θ.
Note that the horizontal lines in b) and c) connect the values of θ -180° and +180°, extremes between which the angle is shown.

From Fig. 5c, it is clear that $X(t)/R(t) = \cos\theta$ e $Y(t)/R(t) = \sin\theta$. So it is possible to write:

$$X(t) = R(t)\cos\theta = R(t)\cos\left(\frac{2\pi t}{T_p}\right)$$
$$Y(t) = R(t)\sin\theta = R(t)\sin\left(\frac{2\pi t}{T_p}\right) \quad \text{eq.2}$$

where $T_p$ is the period and a uniform rotational speed on the $X, Y$ plane has been assumed to replace $\theta$ with $\frac{2\pi t}{T_p}$. This means that we try to reproduce the spatial form of the trajectories, but not necessarily the speed at which they are traveled.
The trend of the radius in each revolution is shown in Fig. 5b, black trace. Its shape suggests the function of the first level of the harmonic oscillator, suitably modified. The expression to reproduce the radius is obtained with the following procedure. First: the function of the first level of the harmonic oscillator (excluding constants) (Fig. 6a), is added to a sigmoid that passes from the $r_i$ value to the $r_{i+1}$ value of two successive radii provided by equation 1. See Fig. 6b at the top for sigmoid and at the bottom for the sum of sigmoid and first level of the harmonic oscillator. Second: to make the function asymmetric as shown in Fig. 5b, the function of level 0 of the harmonic oscillator (excluding constants) is added, slightly shifted from the center point. See Fig. 6c at the top for the function of level 0 of the harmonic oscillator and below for the resulting sum.

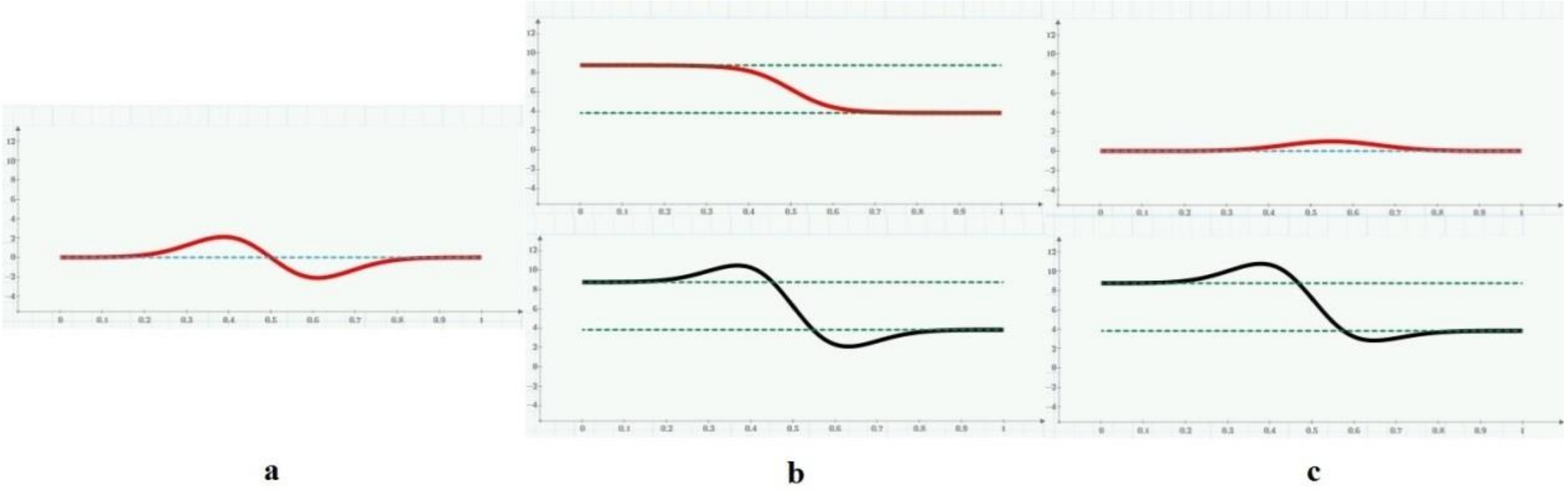

Fig.6. The abscissa are in 2π units. **a**) first level function of the harmonic oscillator:

$-m_5\left[a\left(\frac{t}{T_p}-m_6\right)\right]e^{-\frac{\left[a\left(\frac{t}{T_p}-m_6\right)\right]^2}{2}}$, where $a = 9$, $m_5 = 3.5$, $m_6 = \frac{1}{2}$; the zero line is dashed in blue. **b**) above: sigmoid $m_1 + \frac{m_2}{1+e^{-\frac{\frac{t}{T_p}-m_3}{m_4}}}$ that joins the values $r_i = 8.75$ and $r_{i+1} = 3.828$ in this example, where $m_1 = r_i$, $m_2 = r_{i+1} - r_i$, $m_3 = \frac{1}{2}$, $m_4 = .05$ ; the green dotted lines highlight the values of $r_i$ and $r_{i+1}$. Below is the sum between the function in a) and the sigmoid. **c**) above zero level function of the harmonic oscillator: $m_7 e^{-\frac{\left[a\left(\frac{t}{T_p}-m_8\right)\right]^2}{2}}$, with $m_7 = 1$ and $m_8 = .55$. Below the sum (black trace) of all three functions.

It is thus obtained:

$$m_1 = r_i$$

$$m_2 = r_{i+1} - r_i$$

$$R(t) = m_1 + \frac{m_2}{1+e^{-\frac{\frac{t}{T_p}-m_3}{m_4}}} - m_5\left[a\left(\frac{t}{T_p}-m_6\right)\right]e^{-\frac{\left[a\left(\frac{t}{T_p}-m_6\right)\right]^2}{2}} + m_7 e^{-\frac{\left[a\left(\frac{t}{T_p}-m_8\right)\right]^2}{2}}$$

eq.3

The expression depends on two successive radii, the time $t$ is 0 at $r_i$ and becomes $T_p$ when $r_{i+1}$ is reached. Overall $R(t)$ is a piecewise function being the assembly of expressions, one for each pair of successive radii.

To complete the three-dimensional expression of the curves, we study the trend of the coordinate $Z(t)$ which represents the rise with respect to the $X, Y$ plane (red trace in Fig. 5b). The rise is centered around a θ-value of 35° and concerns about half of the trajectory. For the other half $Z(t)$ is 0, i.e. the trajectory lies on the $X, Y$ plane. Plotting the $Z(t)$ coordinate as a function of the $R(t)$ radius shows that the rise is greater when the highest $R$ values are reached in the trajectory (Fig. 7 reports as an example the case of the Rössler attractor, for $c = 4.00$).
A function with these characteristics is:

$$Z(t) = c_3\left[\frac{r_i}{10}e^{-\left(\frac{2\pi t}{T_p}-\frac{\pi}{6}\right)}\right]^4 \qquad \text{eq. 4}$$

where $c_3 = \frac{1}{3}$ . The phase shift of $-\frac{\pi}{6}$ has the purpose of obtaining the elevation in the position analogous to the Rössler trajectory. The same considerations as in the case of $R(t)$ apply to time $t$. Coefficients are introduced into this and the previous equation to approach the values of the Rössler attractor trajectory.

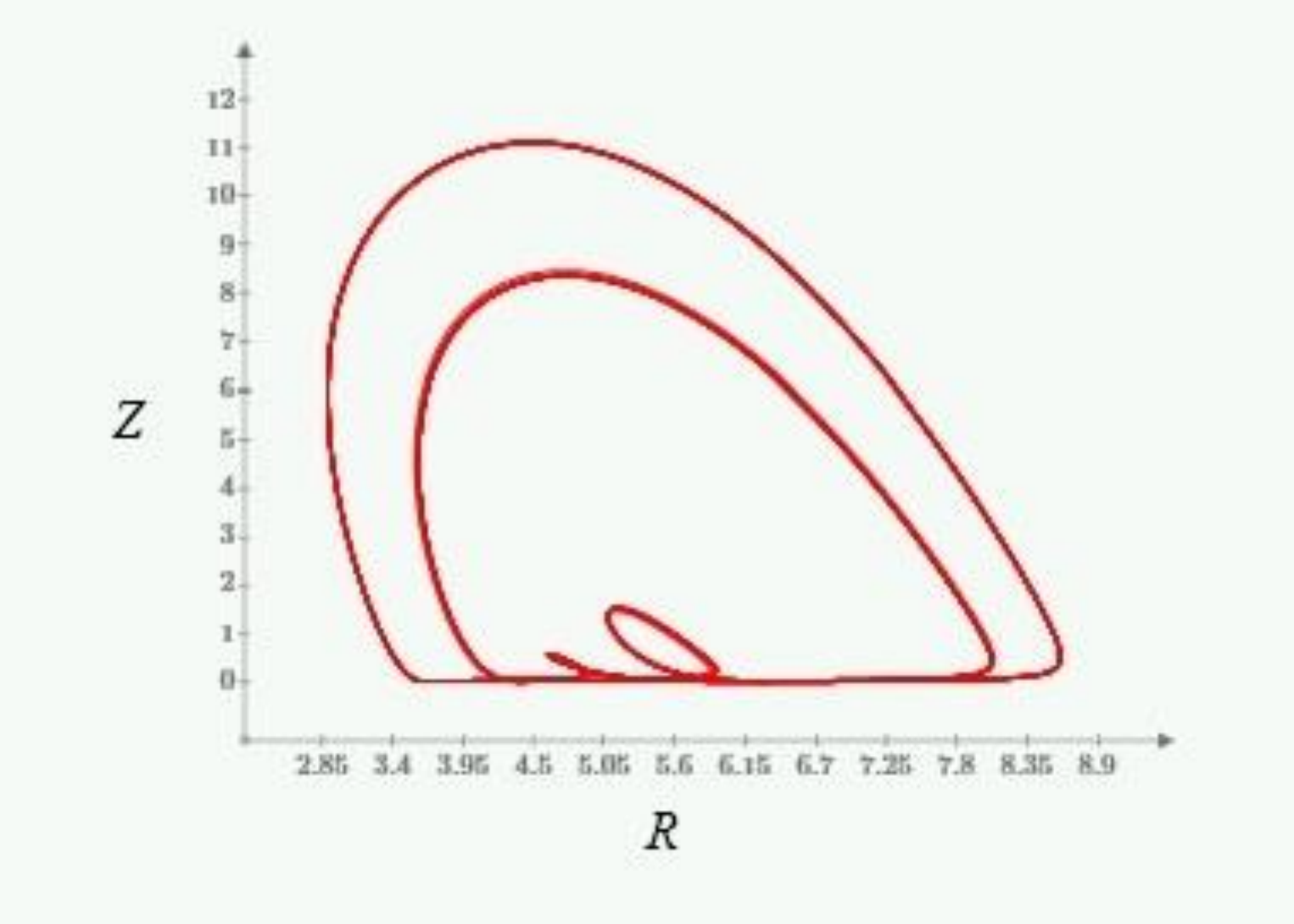

Fig. 7. $Z(t)$ coordinate as a function of radius $R(t)$ in the case of the Rössler attractor, for $c = 4.00$.

In the search for the expressions for the $X$, $Y$, $Z$ components of the trajectories, we focused on expressions with zero derivative at the extremes, in order not to have discontinuity problems in the construction of piecewise functions.

*Results*

Grouping the equations obtained so far:

carrier function

$$x_i = \lambda x_{i-1}(1 - x_{i-1})$$
$$r_i = 10x_i$$

coordinates

$$m_1 = r_i$$
$$m_2 = r_{i+1} - r_i$$
$$R(t) = m_1 + \frac{m_2}{1 + e^{-\frac{\frac{t}{T_p} - m_3}{m_4}}} - m_5\left[a\left(\frac{t}{T_p} - m_6\right)\right]e^{-\frac{\left[a\left(\frac{t}{T_p} - m_6\right)\right]^2}{2}} + m_7 e^{-\frac{\left[a\left(\frac{t}{T_p} - m_8\right)\right]^2}{2}}$$
$$X(t) = -R(t)\cos\left(\frac{2\pi t}{T_p}\right)$$
$$Y(t) = -R(t)\sin\left(\frac{2\pi t}{T_p}\right)$$
$$Z(t) = c_3\left[\frac{r_i}{10} e^{-\left(\frac{2\pi t}{T_p} - \frac{\pi}{6}\right)}\right]^4$$

eq. 5

Note that in the expression of $X(t)$ and $Y(t)$ a π phase has been introduced to change their sign, and thus present the trajectories of the Rössler-like attractor in the same position as the Rössler attractor with respect to the reference system.

$\lambda = 3.5$ gives an Rössler-like attractor of period-4. We report the trend of the radius $R(t)$ on the $X, Y$ plane and of the value $Z(t)$ as a function of θ in Fig. 8.

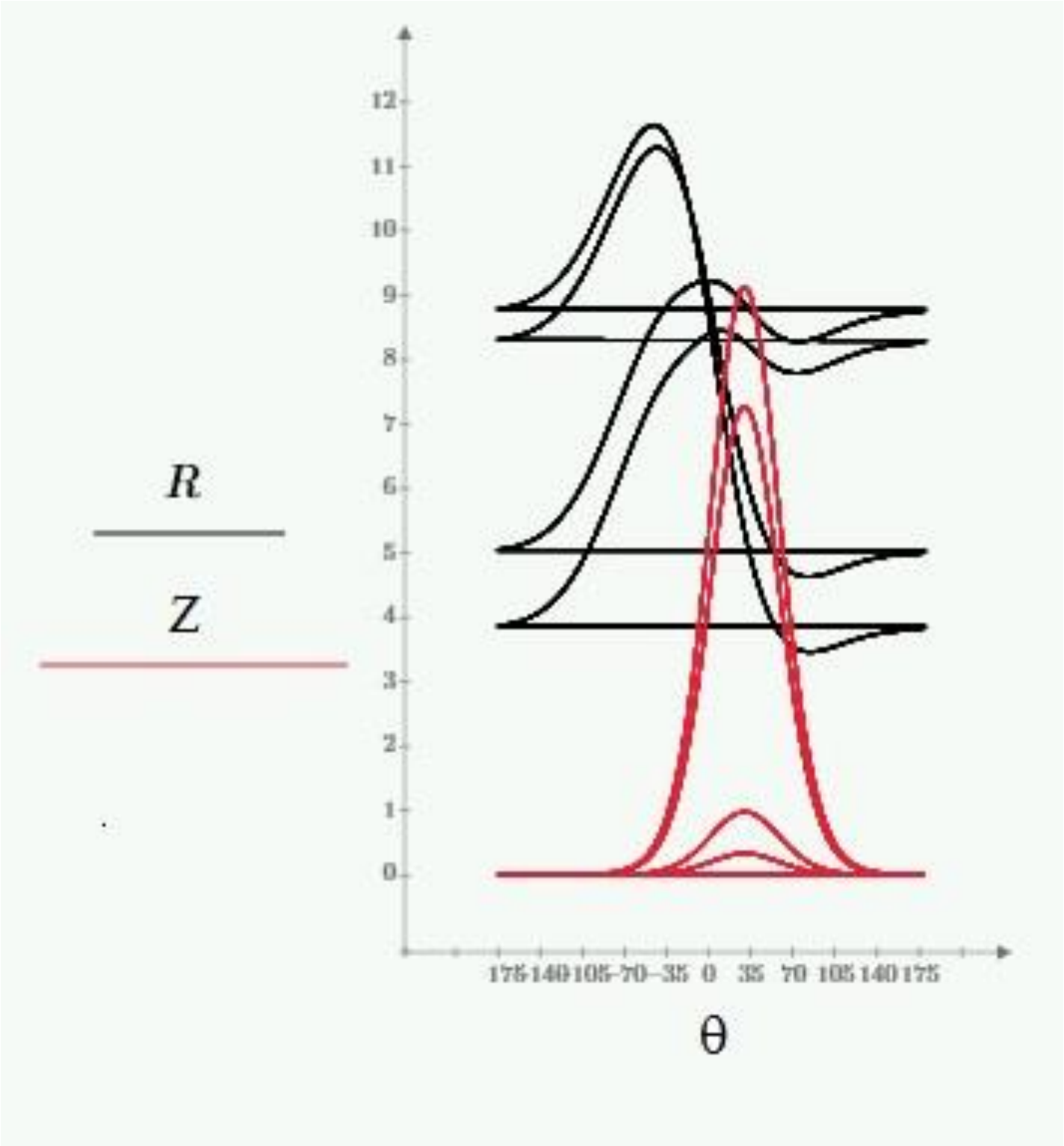

Fig.8 Trend of the radius $R(t)$ on the $X, Y$ plane and of the value $Z(t)$ as a function of θ for the Rössler-like attractor, period-4, for $\lambda = 3.5$.

Note the similarity of the pattern of functions, albeit with different numerical values, with the Rössler attractor of the same periodicity (see Fig. 5b). From this figure it can be graphically inferred that the trajectory does not intersect itself, as it must be for the non-intersection theorem[6]. In fact, there are values of θ for which the radius function $R(t)$ intersects itself, which means that there are pairs of values $X, Y$ equal. But these intersections occur in the region of θ where $Z(t)$ rises from the plane and has no intersections. So there are no triplets of coincident $X$, $Y$, $Z$ values. What has been observed for this particular case is also found in other controlled periodicities, and in the chaotic regime.

Fig. 9 compares the Rössler trajectories with those obtained with the functions described here, for the periodicities 2, 4, 8 and the chaotic regime, obtained by varying the parameter $c$ and the bifurcation parameter $\lambda$ respectively.

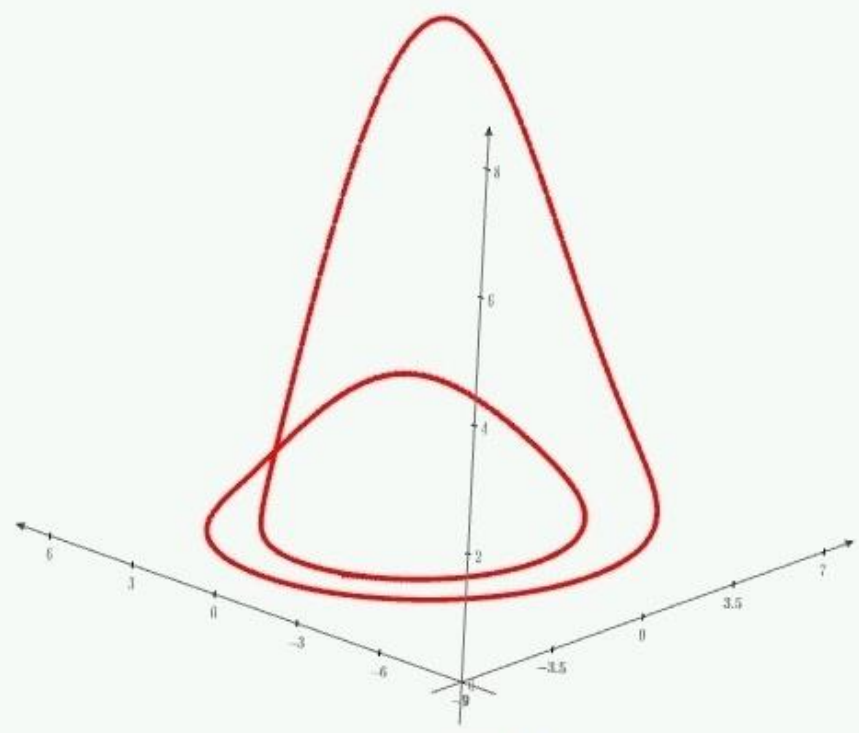

c = 3.25

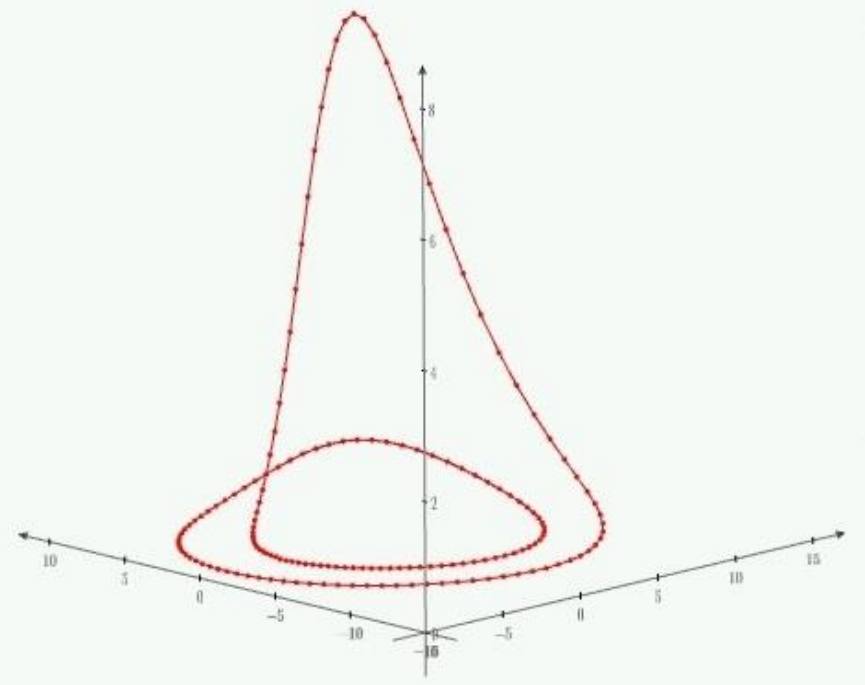

λ = 3.30

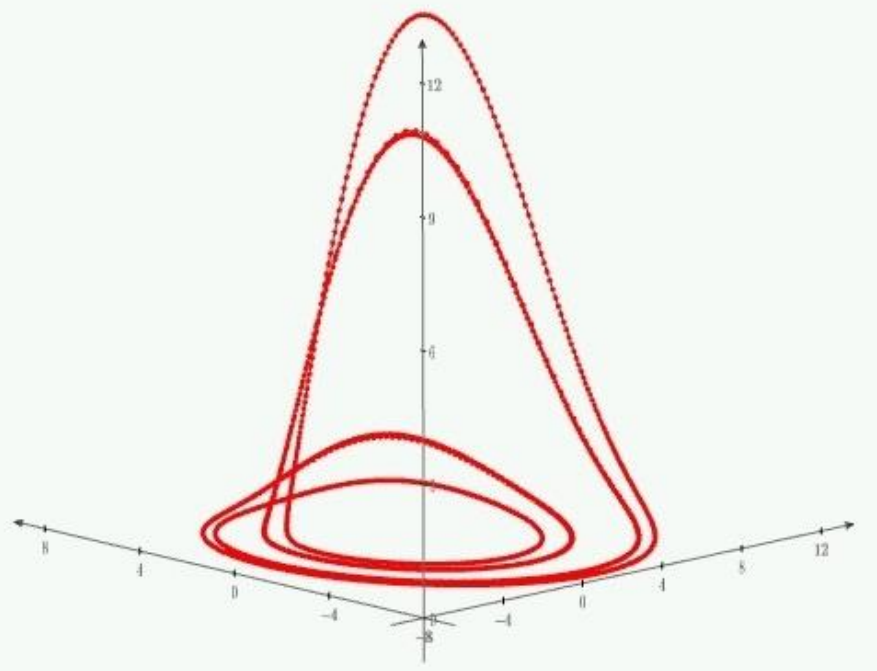

c = 4.00

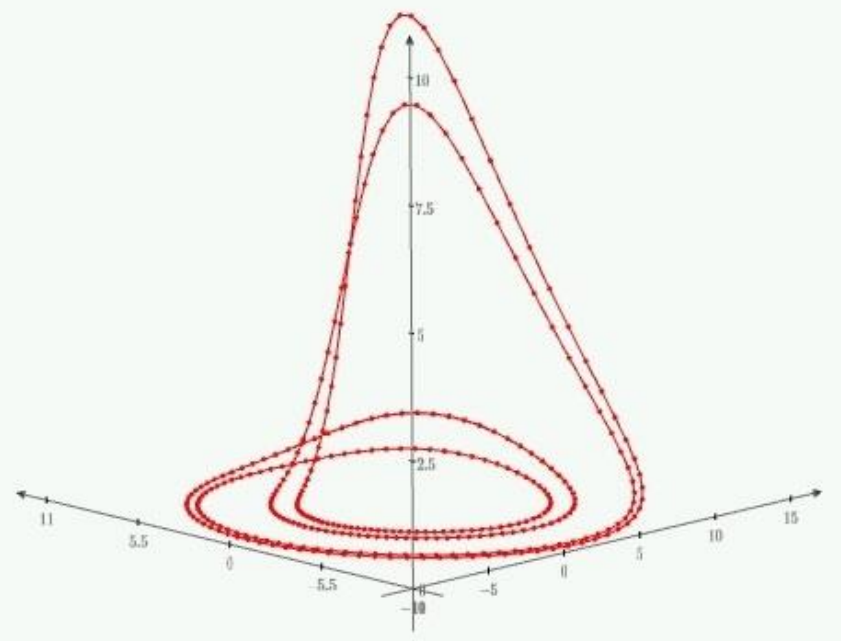

λ = 3.50

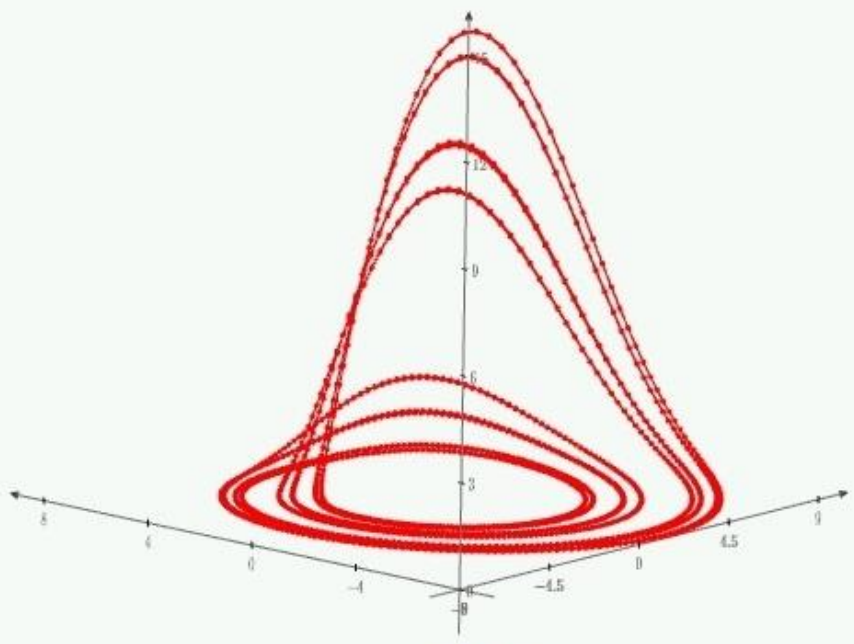

c = 4.20

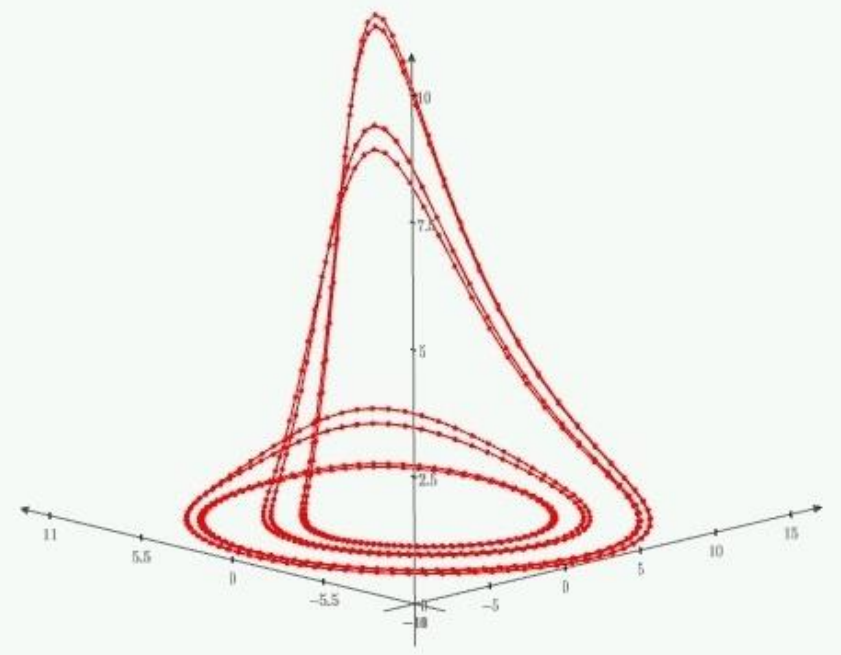

λ = 3.55

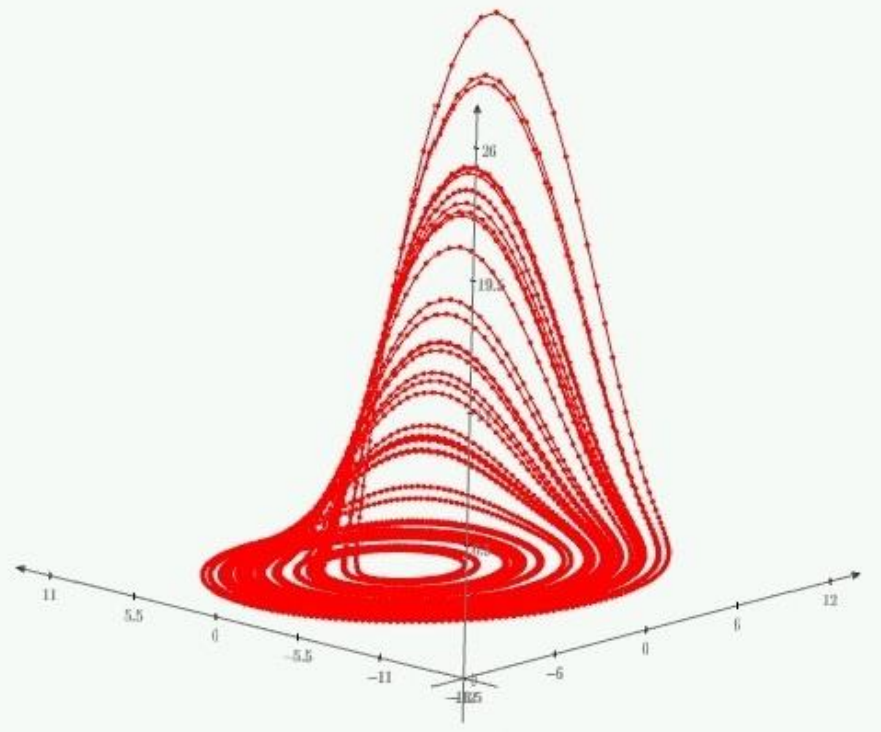

c = 5.70

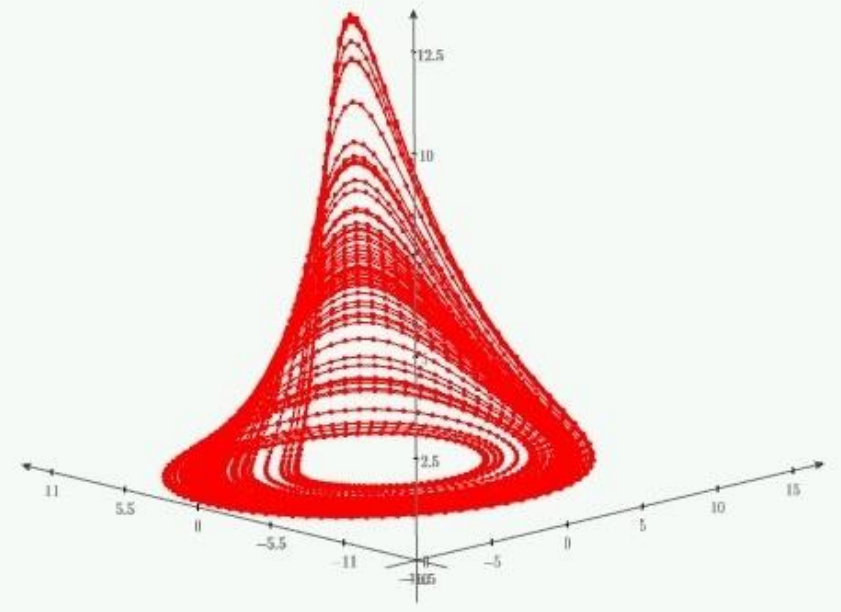

λ = 3.70

Fig.9 From top to bottom, the periodicities 2, 4, 8, and the chaotic regime are shown. Left column: Rössler attractor for the indicated values of parameter $c$. Right column: Rössler-like attractor for the indicated values of $\lambda$.

From the visual inspection of trajectories it seems possible to affirm the topological identity of couples with equal periodicity.

## *Conclusions*

Using iteration of a carrier function and piecewise functions for coordinate components, it was possible to represent Rössler-like trajectories. They are not the analytical expressions of Rössler trajectories, however they are the expression of a nonlinear attractor, with periodic and chaotic regimes. The mathematical approach, different from that hypothesized by Poincaré, allows this realization in apparent contrast to his theorem. This expression differs from the one obtainable with the HAM method, as it is not a series solution.

The following possible developments can be expected.

- using other functions and varying the parameters, it is possible to get closer to the numerical values of the Rössler trajectories.
- trajectories similar to other attractors can be described with the same approach, that is, by setting a carrier function on the basis of the first-return map, and by simulating the trend of the trajectories with appropriate functions.
- non-linear systems whose differential equations are not known, but whose experimental data have been collected can be treated with this approach, since they consist of a succession of experimental values, and therefore conceptually not different from the succession of points obtained by numerical integration of the differential equations. It is possible to build a first-return map, and to study the trend of the trajectories.

*Mathcad Prime 3.0 worksheets*

Programs, written in Mathcad 3.0, with which trajectories were derived. S1: trajectory of the period-4 Rössler attractor, solved numerically. S2 analytical trajectory of the Rössler-like attractor.

S1: *Rössler_c=4_00_program.pdf*

Rossler c=4_00_t=256_ODE

$a := .2$

$b := .2$ $\qquad x := \begin{bmatrix} -6 \\ .16 \\ 0 \end{bmatrix}$

$c := 4.0$

$$D(t,x) := \begin{bmatrix} -x_1 - x_2 \\ x_0 + a \cdot x_1 \\ b + (x_0 - c) \cdot x_2 \end{bmatrix} \qquad Jac(t,x) := \begin{bmatrix} 0 & -1 & -1 \\ 1 & a & 0 \\ x_2 & 0 & x_0 - c \end{bmatrix}$$

$lt := 256$

$lp := 32$

$tol := 10^{-15}$

$$Z := \text{AdamsBDF}(x, 0, lt, lt \cdot lp, D, Jac, tol)$$

$$n := 0 .. lt \cdot lp - 1$$

X, Y projection

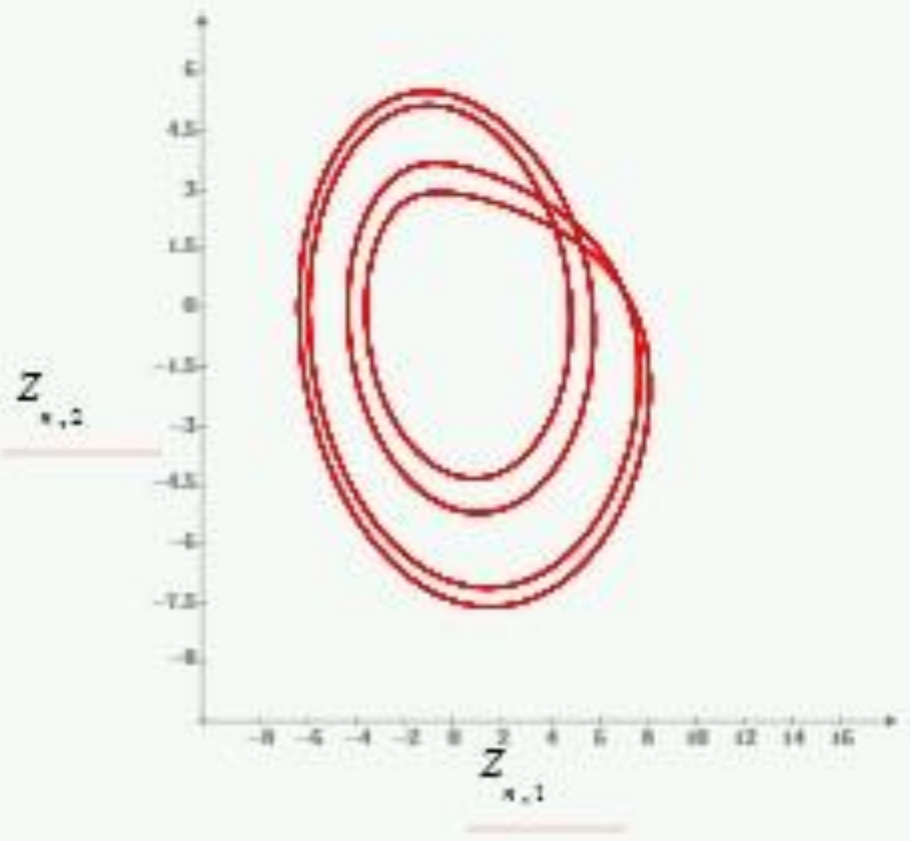

$f1_n := Z_{n,1}$ $\quad$ $f2_n := Z_{n,2}$ $\quad$ $f3_n := Z_{n,3}$

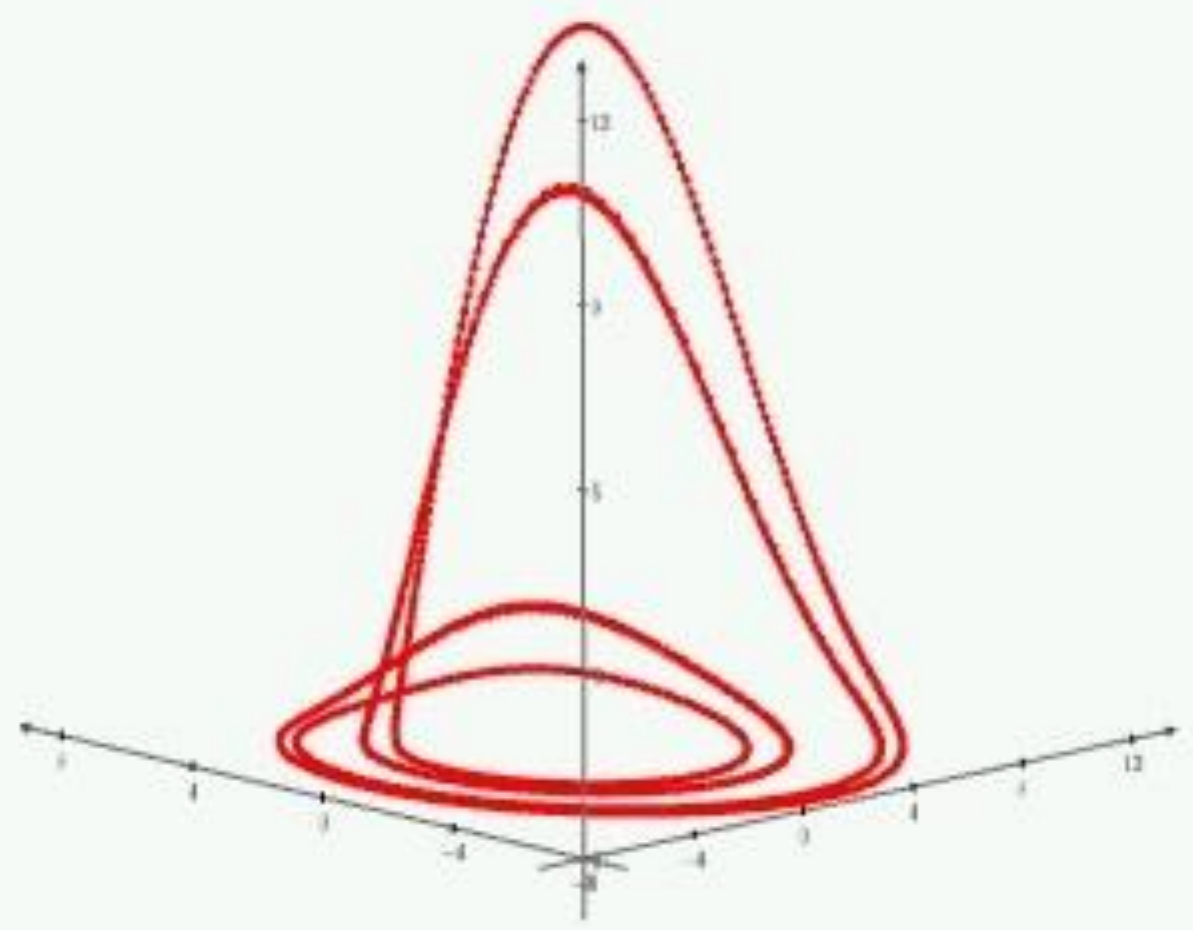

$\begin{bmatrix} f1 \\ f2 \\ f3 \end{bmatrix}$

$$farctg2_n := \frac{\text{atan2}\left(Z_{n,1}, Z_{n,2}\right) \cdot 360}{2 \cdot \pi}$$

$$R_n := \sqrt{{Z_{n,1}}^2 + {Z_{n,2}}^2}$$

Z and R vs θ

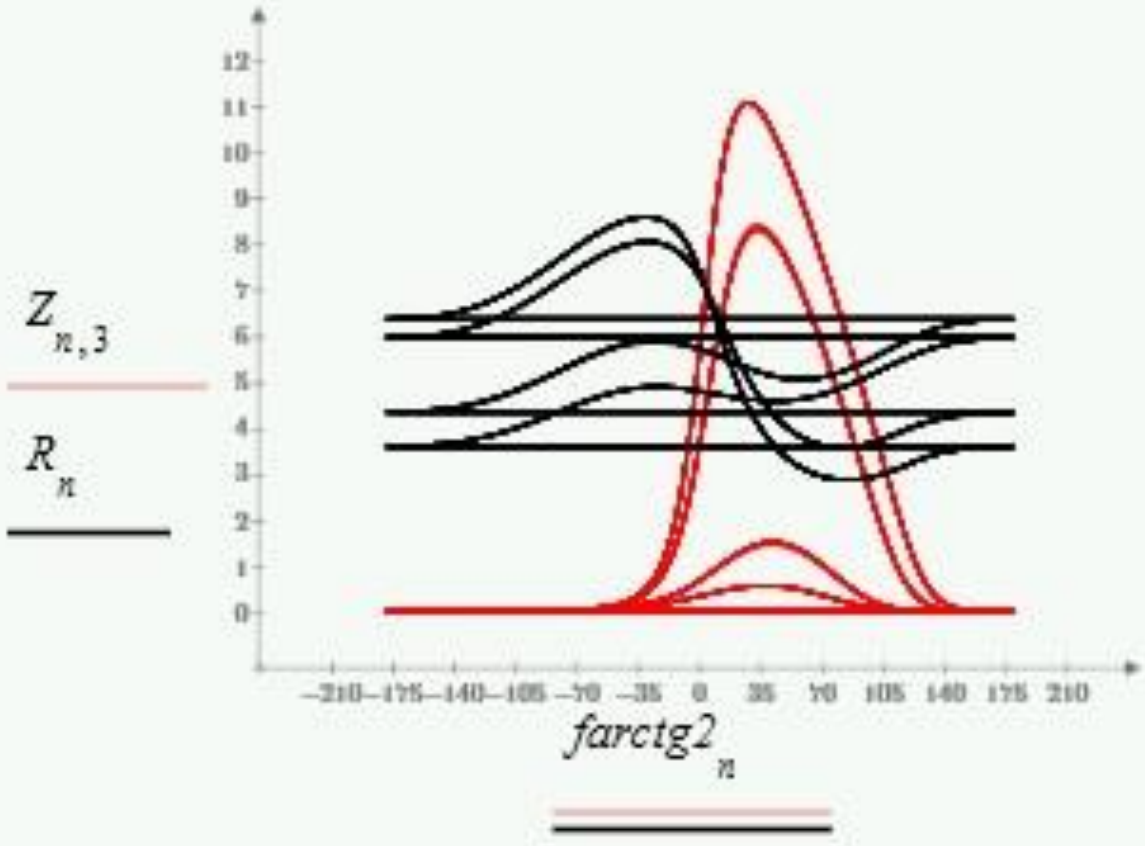

X, Y, X/R, Y/R vs θ

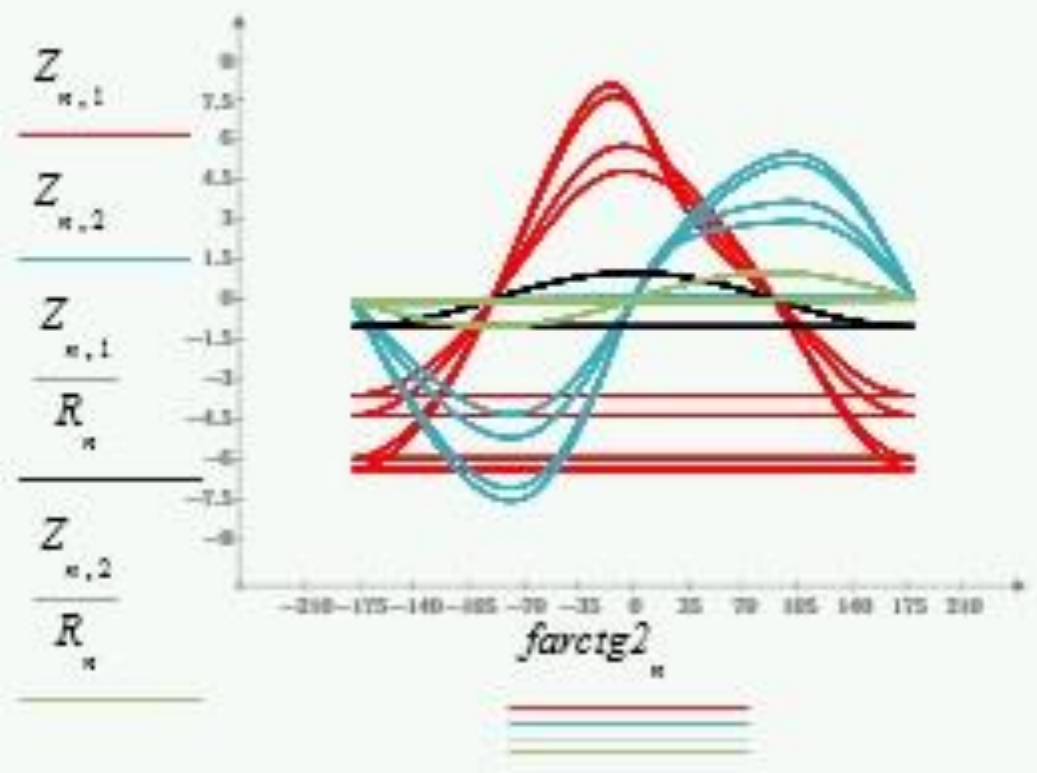

Z vs R

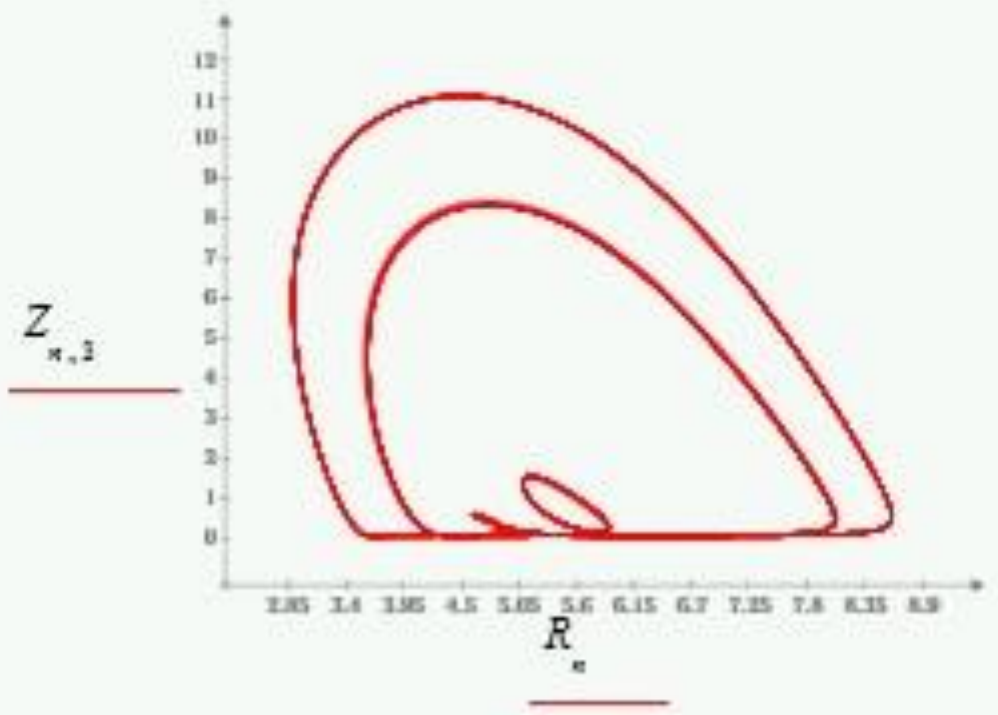

S2: *Rössler-like_lambda=3_50_program.pdf*

assembledRossler_osc_lambda=3_5

bifurcation parameter $\lambda$

niter=number of iterations

$\lambda := 3.5$

$niter := 64$

iteration of the logistic map; initial value xiter=.5

npoints is the number of point for cicle; Tp is the period

$npoints := 80$ $\quad Tp := npoints$

$Tp = 80$

iterlog is a vector containing the iteration values of the logistic map

$$iterlog := \left\| \begin{array}{l} xiter \leftarrow .5 \\ \text{for } i \in 0..niter \\ \quad \left\| \begin{array}{l} iter_log_i \leftarrow xiter \cdot 10 \\ xiter \leftarrow \lambda \cdot xiter \cdot (1 - xiter) \end{array} \right. \\ iter_log \end{array} \right. = \begin{bmatrix} 5 \\ 8.75 \\ 3.828 \\ 8.269 \\ 5.009 \\ 8.75 \\ 3.828 \\ 8.269 \\ 5.009 \\ \vdots \end{bmatrix}$$

connection between consecutive radii obtained by first osc + logistic+zero osc functions

$t := 0..npoints$

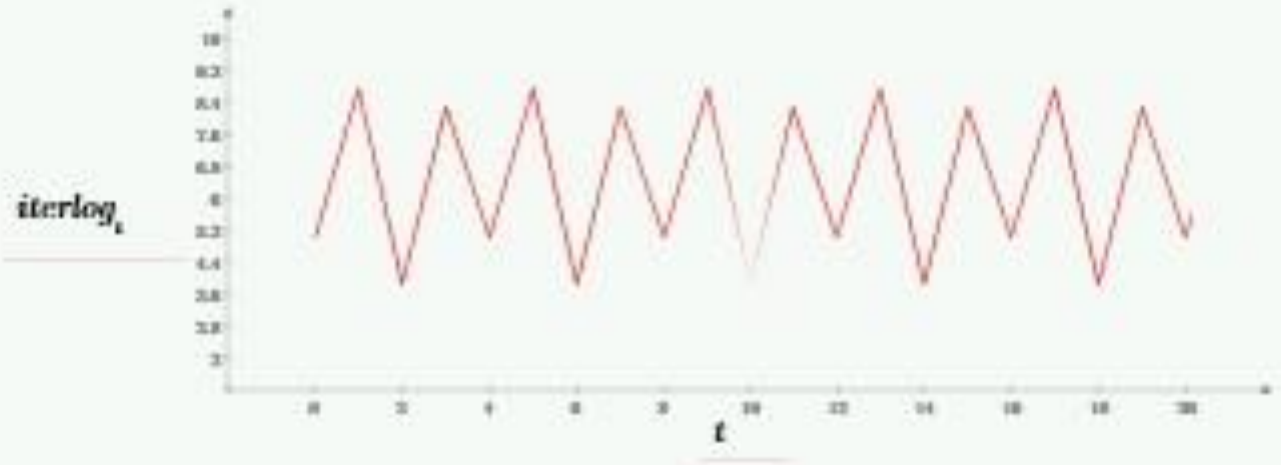

nshare: number of assembled pieces; nshare max=niter-1; coord: matrix of the coordinates

internal j is introduced so to have integer indexes

$$nshare := niter - 1$$

$$m3 := .5 \qquad m4 := .05 \quad m5 := 3.5 \qquad m6 := .5 \quad m7 := 2.5 \quad m8 := .55$$

$$a := 7 \qquad c3 := \frac{1}{3.5}$$

$$
\begin{aligned}
coord := \; & j \leftarrow -1 \\
& \text{for } i \in 0 \,..\, nshare \\
& \quad \text{for } t \in 0 \,..\, npoints - 1 \\
& \quad\quad j \leftarrow j + 1 \\
& \quad\quad m1 \leftarrow iterlog(i) \\
& \quad\quad m2 \leftarrow iterlog(i+1) - iterlog(i) \\
& \quad\quad logistic_t \leftarrow m1 + \frac{m2}{\left(1 + \exp\left(-\frac{\left(\frac{t}{Tp} - m3\right)}{m4}\right)\right)} \\
& \quad\quad osc_t \leftarrow -m5 \cdot \left(a \cdot \left(\frac{t}{Tp} - m6\right)\right) \cdot \exp\left(-\frac{\left(a \cdot \left(\frac{t}{Tp} - m6\right)\right)^2}{2}\right) \\
& \quad\quad osc0_t \leftarrow m7 \cdot \exp\left(-\frac{\left(a \cdot \left(\frac{t}{Tp} - m8\right)\right)^2}{2}\right) \\
& \quad\quad radius_t \leftarrow logistic_t + osc_t + osc0_t \\
& \quad\quad coord_int_{j,0} \leftarrow -\cos\left(\frac{2 \cdot \pi \cdot t}{Tp}\right) \cdot radius_t \\
& \quad\quad coord_int_{j,1} \leftarrow -\sin\left(\frac{2 \cdot \pi \cdot t}{Tp}\right) \cdot radius_t \\
& \quad\quad coeff_t \leftarrow \frac{iterlog(i)}{10} \\
& \quad\quad coord_int_{j,2} \leftarrow c3 \cdot \left(coeff_t \cdot \exp\left(-\cos\left(\frac{2 \cdot \pi \cdot t}{Tp} - \frac{2\,\pi}{12}\right)\right)\right)^4 \\
& \quad coord_int \\
& coord_int
\end{aligned}
$$

$$coord = \begin{bmatrix} -5.028 & 0 & 5.589\cdot10^{-4} \\ -5.022 & -0.395 & 4.829\cdot10^{-4} \\ -4.987 & -0.79 & 4.266\cdot10^{-4} \\ -4.924 & -1.182 & 3.856\cdot10^{-4} \\ -4.834 & -1.571 & 3.569\cdot10^{-4} \\ -4.718 & -1.954 & 3.385\cdot10^{-4} \\ -4.576 & -2.331 & 3.289\cdot10^{-4} \\ -4.409 & -2.702 & 3.275\cdot10^{-4} \\ -4.218 & -3.064 & 3.343\cdot10^{-4} \\ -4.003 & -3.419 & 3.497\cdot10^{-4} \\ -3.765 & -3.765 & 3.748\cdot10^{-4} \\ -3.505 & -4.103 & 4.114\cdot10^{-4} \\ & & \vdots \end{bmatrix}$$

$j := \mathrm{rows}(coord)$ $\qquad j = 5.12\cdot10^{3}$

ttot: time of the reconstructed trajectory

$ttot := 0..j-1$

$X1 := coord^{\langle 0\rangle}$ $Y1 := coord^{\langle 1\rangle}$ $\qquad Z1 := coord^{\langle 2\rangle}$

x, y projection

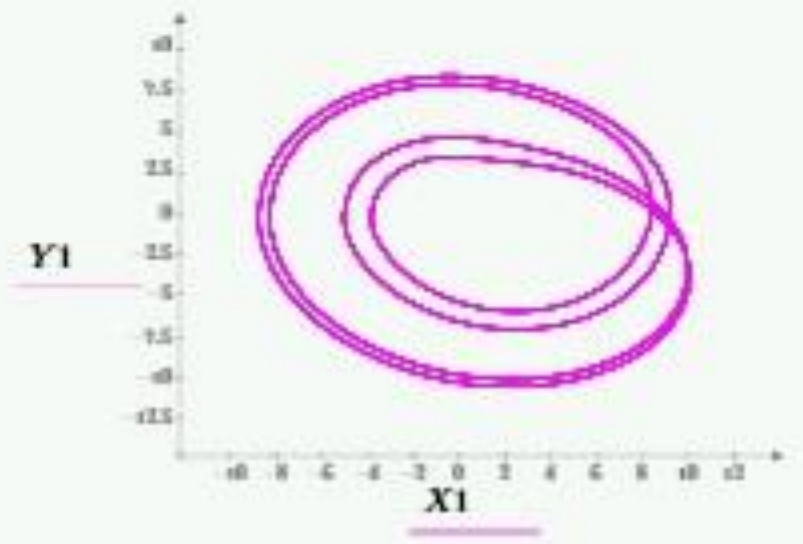

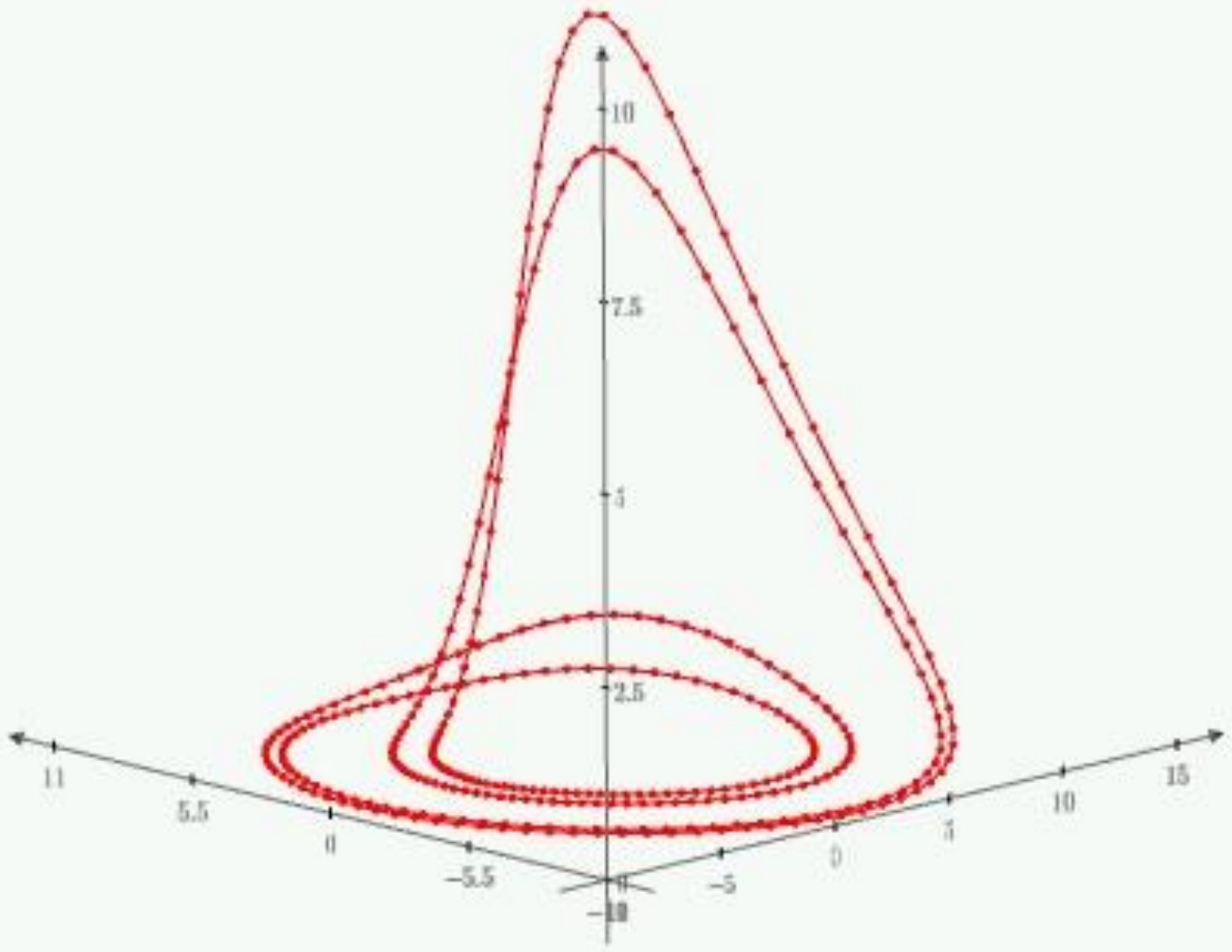

$$\begin{bmatrix} X1 \\ Y1 \\ Z1 \end{bmatrix}$$

$$farctg2_{ttot} := \frac{\text{atan2}\left(coord_{ttot,0}, coord_{ttot,1}\right) \cdot 360}{2 \cdot \pi}$$

$$R_{ttot} := \sqrt{coord_{ttot,0}^{\ 2} + coord_{ttot,1}^{\ 2}}$$

Z vs θ

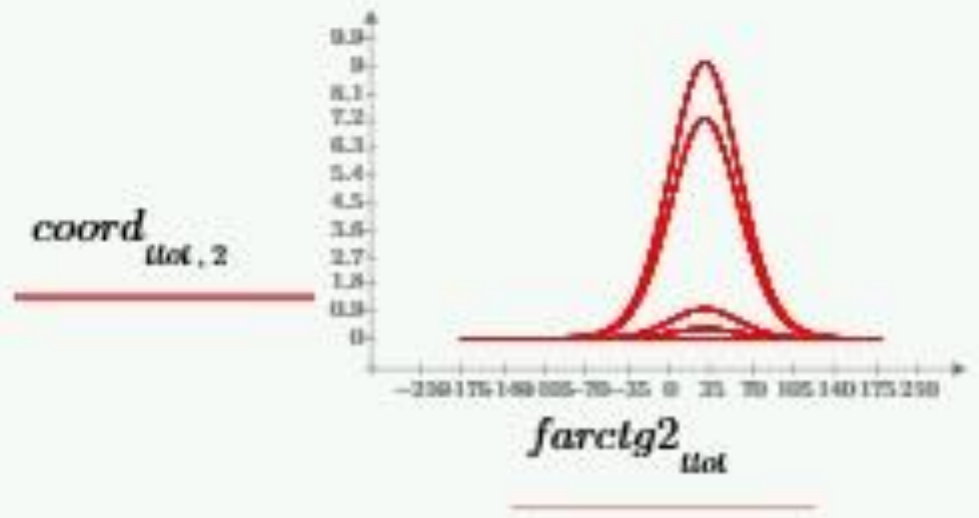

R vs θ

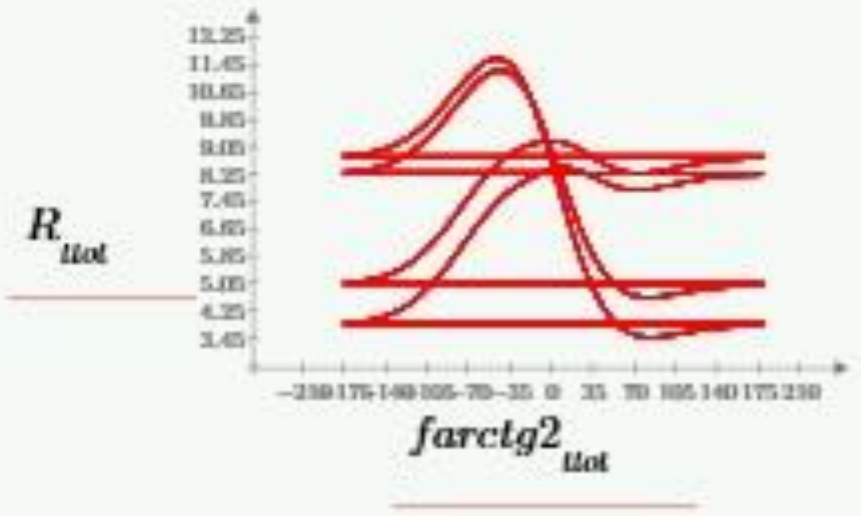

R and Z vs θ

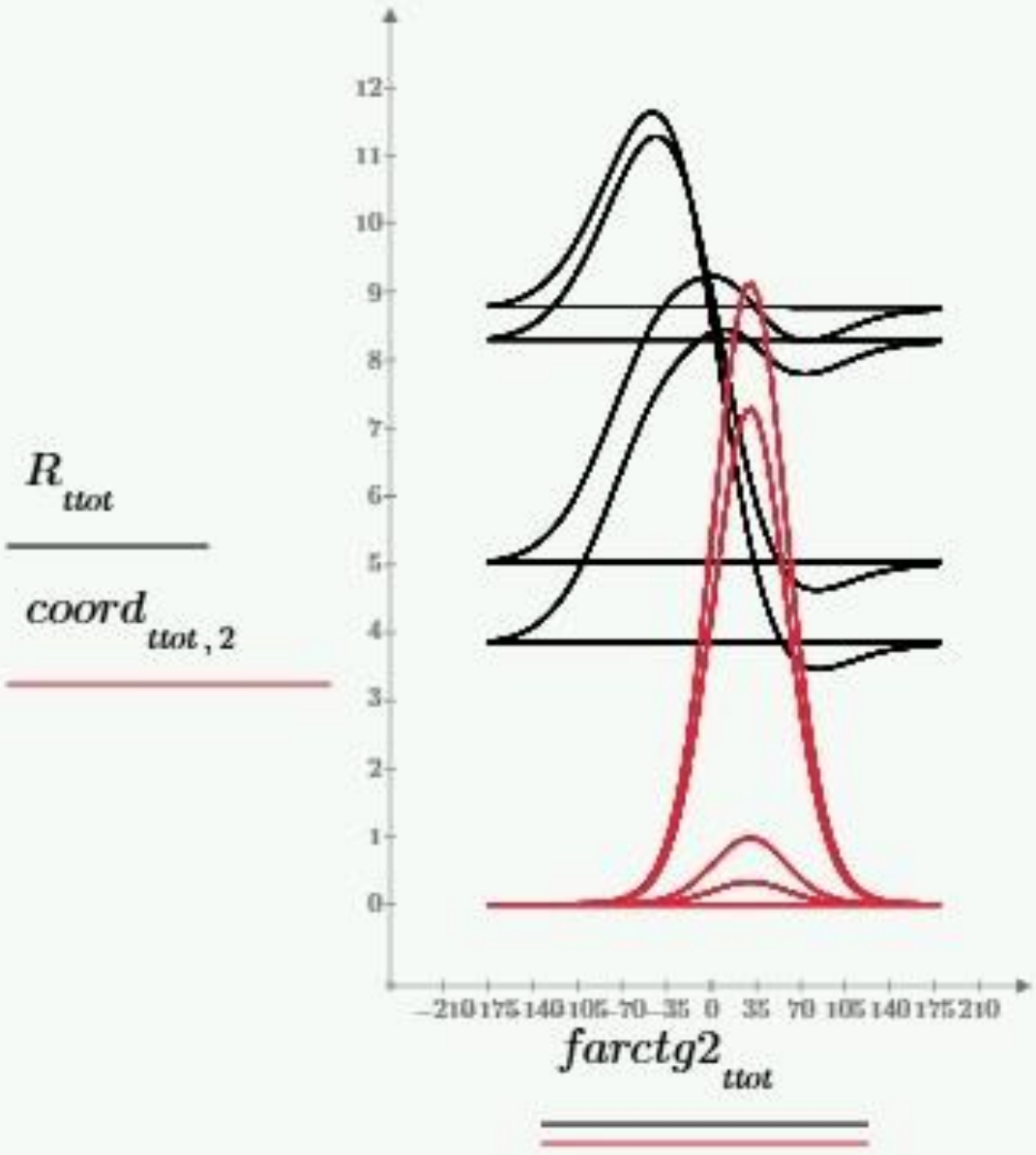